\documentclass[12pt, bibalpha]{asl}

\usepackage{amssymb}
\usepackage{graphicx}
\usepackage{hyperref}
\usepackage{ifthen}
\usepackage[bbgreekl]{mathbbol}
\usepackage{mathrsfs}
\usepackage{pict2e}
\usepackage{xargs}
\usepackage{xspace}

\title{Recurrence and the existence of invariant measures}



\author{Manuel J. Inselmann}
\revauthor{Inselmann, Manuel J.}
\address{
  Manuel Inselmann \\
  Kurt G\"{o}del Research Center for Mathematical Logic \\
  Universit\"{a}t Wien \\
  W\"{a}hringer Stra{\ss}e 25 \\
  1090 Wien \\
  Austria
}
\email{manuel.inselmann@univie.ac.at}

\author{Benjamin D. Miller}
\revauthor{Miller, Benjamin D.}
\address{
  Benjamin D. Miller \\
  Kurt G\"{o}del Research Center for Mathematical Logic \\
  Universit\"{a}t Wien \\
  W\"{a}hringer Stra{\ss}e 25 \\
  1090 Wien \\
  Austria
 }
\email{benjamin.miller@univie.ac.at}
\urladdr{\url{http://www.logic.univie.ac.at/benjamin.miller}}

\thanks{The authors were supported in part by FWF Grants
  P28153 and P29999.}

\DeclareSymbolFontAlphabet{\amsmathbb}{AMSb}

\newcommand{\definedterm}[1]{\emph{#1}}

\newcommand{\action}{\curvearrowright}
\newcommand{\Bairespace}[1][]{
  \ifthenelse{\equal{#1}{}}{\functions{\N}{\N}}{\functions{#1}{\N}}
}
\newcommand{\Bairetree}[1][]{
  \ifthenelse{\equal{#1}{}}{\functions{<\N}{\N}}{\functions{#1}{\N}}
}
\newcommand{\ball}[3][]{\calB_{#1}(#2, #3)}
\newcommand{\bbg}{\mathbb{g}}
\newcommand{\bbk}{\mathbb{k}}
\newcommand{\bbs}{\mathbb{s}}

\newcommand{\calB}{\mathcal{B}}
\newcommand{\calF}{\mathcal{F}}

\newcommand{\calN}{\mathcal{N}}
\newcommand{\calO}{\mathcal{O}}
\newcommand{\calP}{\mathcal{P}}
\newcommand{\calS}{\mathcal{S}}
\newcommand{\calT}{\mathcal{T}}

\newcommand{\Cantorspace}[1][]{
  \ifthenelse{\equal{#1}{}}{\functions{\N}{2}}{\functions{#1}{2}}
}
\newcommand{\Cantortree}[1][]{
  \ifthenelse{\equal{#1}{}}{\functions{<\N}{2}}{\functions{#1}{2}}
}

\newcommand{\closedsets}[1]{F(#1)}
\newcommand{\closure}[1]{\overline{#1}}

\newcommand{\composition}{\circ}

\newcommandx{\concatenation}[2][1 = undefined, 2 = undefined]{
  \ifthenelse{\equal{#1}{undefined}}{{}\smallfrown}{
    \ifthenelse{\equal{#2}{undefined}}{\bigoplus #1}{\bigoplus_{#1} #2}
  }
}
\newcommand{\constantsequence}[2]{\sequence{#1}^{#2}}
\newcommandx{\convolution}[2][1 = undefined, 2 = undefined]{
  \ifthenelse{\equal{#1}{undefined}}{\mathrel{*}}{
    \ifthenelse{\equal{#2}{undefined}}{\bigotimes #1}{\bigotimes_{#1} #2}
  }
}

\newcommandx{\Deltaclass}[2][1=,2=]{
  \ifthenelse{\equal{#2}{}}{\mathbf{\Delta}_{#1}}{\mathbf{\Delta}^{#1}_{#2}}
}
\newcommand{\dzero}[1][]{\ifthenelse{\equal{#1}{}}{\mathbb{d}_0^{\Delta}}{\mathbb{d}_0^{\Delta}(#1)}}

\newcommandx{\disjointunion}[2][1 =, 2 =]{
  \ifthenelse{\equal{#1}{}}{\sqcup}{
    \ifthenelse{\equal{#2}{}}{\bigsqcup #1}{{\bigsqcup_{#1} #2}}
  }
}

\newcommand{\Dtwo}[1]{#1 \setminus #1}

\newcommand{\dual}[1]{\breve{#1}}

\newcommand{\equivalenceclass}[2]{[#1]_{#2}}

\newcommand{\extensions}[1]{\calN_{#1}}
\newcommand{\Ezero}[1][]{
  \ifthenelse{\equal{#1}{}}{\amsmathbb{E}_0}{\amsmathbb{E}_0(#1)}
}

\newcommand{\forcomeagerlymany}{\forall^*}
\newcommand{\from}{\colon}

\newcommandx{\functions}[3][3 =]{
  \ifthenelse{\equal{#3}{}}{#2^{#1}}{#2^{#1}_{#3}}
}
\newcommandx{\gzero}[2][2 =]{\ifthenelse{\equal{#2}{}}{\bbg_0^{#1}}{\bbg_0^{#1}(#2)}}
\newcommand{\gammazero}[1][]{\ifthenelse{\equal{#1}{}}{\bbgamma_0}{\bbgamma_0(#1)}}
\newcommand{\Gdelta}{$G_\delta$\xspace}
\newcommand{\goesto}{\rightarrow}
\newcommand{\graph}[1]{\mathrm{graph}(#1)}

\newcommand{\GzeroN}[1][]{
  \ifthenelse{\equal{#1}{}}{\amsmathbb{G}_0^\N}{\amsmathbb{G}_{0, #1}^{\N}}
}
\newcommand{\hittingtimes}[2]{\Delta(#1, #2)}
\newcommand{\horizontalsection}[2]{#1^{#2}}
\newcommand{\id}{\mathrm{id}}
\newcommand{\identity}[1]{1_{#1}}
\newcommand{\image}[2]{#1(#2)}
\renewcommand{\implies}{\Longrightarrow}
\renewcommand{\impliedby}{\Longleftarrow}

\newcommand{\infimum}[2][]{
  \ifthenelse{\equal{#1}{}}{\inf #1}{\inf_{#1}{#2}}
}
\newcommandx{\intersection}[2][1 =, 2 =]{
  \ifthenelse{\equal{#1}{}}{\cap}{
    \ifthenelse{\equal{#2}{}}{\bigcap #1}{{\bigcap_{#1} #2}}
  }
}
\newcommand{\inverse}[1]{#1^{-1}}

\newcommand{\kzero}[1][]{\ifthenelse{\equal{#1}{}}{\bbk_0}{\bbk_0(#1)}}
\newcommand{\length}[1]{|#1|}
\newcommandx{\limit}[2][1 =, 2 =]{
  \ifthenelse{\equal{#1}{}}{\lim}{
    \ifthenelse{\equal{#2}{}}{\lim #1}{{\lim_{#1} #2}}
  }
}

\newcommand{\mathor}{\text{ or }}

\newcommand{\minimaldecomposition}[2]{F_{#1}^{#2}}

\newcommand{\muzero}{m}
\newcommand{\N}{\amsmathbb{N}}
\newcommand{\odometer}{\sigma}
\newcommand{\openinterval}[2]{(#1, #2)}
\newcommand{\orbit}[2]{#2 #1}

\newcommand{\orbitequivalencerelation}[2]{E_{#1}^{#2}}

\newcommand{\pair}[2]{(#1, #2)}
\newcommandx{\Piclass}[2][1=,2=]{
  \ifthenelse{\equal{#2}{}}{\mathbf{\Pi}_{#1}}{\mathbf{\Pi}^{#1}_{#2}}
}
\newcommand{\powerset}[1]{\calP(#1)}
\newcommand{\preimage}[2]{#1^{-1}(#2)}
\newcommand{\principleF}[1]{\calF_{#1}}

\newcommand{\probabilitymeasures}[1]{P(#1)}
\newcommandx{\product}[2][1 =, 2 =]{
  \ifthenelse{\equal{#1}{}}{\times}{
    \ifthenelse{\equal{#2}{}}{\prod #1}{{\prod_{#1} #2}}
  }
}
\newcommandx{\projection}[2][1 =, 2 =]{
  \ifthenelse{\equal{#1}{}}{\mathrm{proj}}{
    \ifthenelse{\equal{#2}{}}{\projection_{#1}}{
      \image{\projection[#1]}{#2}
    }
  }
}
\newcommand{\pushforward}[2]{#1_* #2}
\renewcommand{\restriction}[2]{#1 \upharpoonright #2}
\newcommand{\returntimes}[2][]{\Delta(#2)}

\newcommand{\rhozero}[2][]{\ifthenelse{\equal{#1}{}}{\bbrho_0^{#2}}
  {\bbrho_0^{#1, #2}}}
\newcommand{\Rzero}[1][]{
  \ifthenelse{\equal{#1}{}}{\amsmathbb{R}_0}{\amsmathbb{R}_0(#1)}
}
\newcommand{\saturation}[2]{[#1]_{#2}}

\newcommandx{\sequence}[2][2 = undefined]{
  \ifthenelse{\equal{#2}{undefined}}{(#1)}{
    (#1)_{#2}
  }
}
\newcommandx{\set}[2][2 = undefined]{
  \ifthenelse{\equal{#2}{undefined}}{\{ #1 \}}{
    \{ #1 \suchthat #2 \}
  }
}
\newcommand{\setcomplement}[1]{\mathord{\sim} #1}
\newcommandx{\sets}[3][3 =]{
  \ifthenelse{\equal{#3}{}}{[#2]^{#1}}{[#2]^{#1}_{#3}}
}

\newcommandx{\Sigmaclass}[2][1=,2=]{
  \ifthenelse{\equal{#2}{}}{\mathbf{\Sigma}_{#1}}{\mathbf{\Sigma}^{#1}_{#2}}
}

\newcommand{\suchthat}{\mid}

\newcommand{\symmetricdifference}{\mathrel{\triangle}}
\newcommand{\szero}[1][]{\ifthenelse{\equal{#1}{}}{\bbs_0}{\bbs_0(#1)}}

\newcommandx{\union}[2][1 =, 2 =]{
  \ifthenelse{\equal{#1}{}}{\cup}{
    \ifthenelse{\equal{#2}{}}{\bigcup #1}{{\bigcup_{#1} #2}}
  }
}
\newcommand{\unions}[1]{#1_\sigma}
\newcommand{\verticalsection}[2]{#1_{#2}}

\newcommand{\Z}{\amsmathbb{Z}}

\newcommand{\Baire}{Baire\xspace}
\newcommand{\Becker}{Beck\-er\xspace}
\newcommand{\Borel}{Bor\-el\xspace}
\newcommand{\Effros}{Eff\-ros\xspace}
\newcommand{\Eigen}{Eigen\xspace}

\newcommand{\Fell}{Fell\xspace}
\newcommand{\Glasner}{Glas\-ner\xspace}
\newcommand{\Glimm}{Glimm\xspace}

\newcommand{\Hajian}{Haj\-i\-an\xspace}
\newcommand{\Harrington}{Har\-ring\-ton\xspace}
\newcommand{\Hausdorff}{Haus\-dorff\xspace}

\newcommand{\Kechris}{Kech\-ris\xspace}
\newcommand{\Kuratowski}{Kur\-at\-ow\-ski\xspace}
\newcommand{\Lebesgue}{Leb\-es\-gue\xspace}

\newcommand{\Louveau}{Lou\-veau\xspace}
\newcommand{\Lusin}{Lu\-sin\xspace}
\newcommand{\Miller}{Mill\-er\xspace}
\newcommand{\Montgomery}{Mont\-gom\-ery\xspace}

\newcommand{\Nadkarni}{Nad\-kar\-ni\xspace}
\newcommand{\Nikodym}{Nik\-o\-d\'ym\xspace}
\newcommand{\Novikov}{No\-vik\-ov\xspace}
\newcommand{\Pettis}{Pet\-tis\xspace}
\newcommand{\Polish}{Po\-lish\xspace}
\newcommand{\Radon}{Ra\-don\xspace}

\newcommand{\Sierpinski}{Sier\-pi\'{n}\-ski\xspace}

\newcommand{\Tserunyan}{Tser\-un\-yan\xspace}
\newcommand{\Ulam}{U\-lam\xspace}

\newcommand{\Weiss}{Weiss\xspace}

\renewcommand{\Box}{
  \begin{picture}(6,6)
    \put(0,0){\line(1,0){6}}
    \put(0,0){\line(0,1){6}}
    \put(6,0){\line(0,1){6}}
    \put(0,6){\line(1,0){6}}
    \put(0,0){\line(1,1){6}}
    \put(6,0){\line(-1,1){6}}
  \end{picture}
}

\newenvironment{lemmaproof}{
  \renewcommand{\qedsymbol}{$\BoxBox$}
  \begin{proof}
}{\end{proof}}

\newenvironment{propositionproof}{
  \renewcommand{\qedsymbol}{$\Box$}
  \begin{proof}
}{\end{proof}}

\newenvironment{theoremproof}{
  \renewcommand{\qedsymbol}{$\Box$}
  \begin{proof}
}{\end{proof}}

\newtheorem{lemma}{Lemma}[section]
\newtheorem{proposition}[lemma]{Proposition}

\newtheorem{theorem}[lemma]{Theorem}

\theoremstyle{definition}

\begin{document}

\begin{abstract}
  We show that recurrence conditions do not yield invariant \Borel
  probability measures in the descriptive set-theoretic milieu, in the
  strong sense that if a \Borel action of a locally compact \Polish
  group on a standard \Borel space satisfies such a condition but
  does not have an orbit supporting an invariant \Borel probability
  measure, then there is an invariant \Borel set on which the action
  satisfies the condition but does not have an invariant \Borel
  probability measure.
\end{abstract}

\maketitle

Suppose that $X$ is a \Borel space and $T \from X \to X$ is a
\Borel automorphism. Given a set $S \subseteq \Z$, a set $Y
\subseteq X$ is \definedterm{$S$-wandering} if $\image{T^m}{Y}
\intersection \image{T^n}{Y} = \emptyset$ for all distinct $m, n \in
S$, and \definedterm{weakly wandering} if there is an infinite set $S
\subseteq \Z$ for which it is $S$-wandering. We say that a set $Y
\subseteq X$ is \definedterm{$T$-complete} if it intersects each
orbit of $T$, and a \Borel probability measure $\mu$ on $X$ is
\definedterm{$T$-invariant} if $\mu(B) = \mu(\image{T}{B})$ for all
\Borel sets $B \subseteq X$.

It is well-known that if $\nu$ is a \Borel probability measure on $X$,
then the inexistence of a weakly-wandering $\nu$-positive set yields
a $T$-invariant \Borel probability measure $\mu \gg \nu$ (see \cite
{Zak93} for the generalization to groups of \Borel
automorphisms). In light of this and his own characterization of the
existence of invariant \Borel probability measures (see \cite
{Nad90}), \Nadkarni asked whether the inexistence of a
weakly-wandering $T$-complete \Borel set yields a $T$-invariant
\Borel probability measure in the special case that $X$ is a standard
\Borel space.

We say that a set $Y \subseteq X$ is \definedterm{$T$-invariant} if
$Y = \image{T}{Y}$, and a \Borel probability measure $\mu$ on $X$
is \definedterm{$T$-ergodic} if every $T$-invariant \Borel set is
$\mu$-conull or $\mu$-null. When $X$ is a standard \Borel space, a
\Borel equivalence relation $F$ on $X$ is \definedterm{smooth} if
there is a standard \Borel space $Z$ for which there is a \Borel
function $\pi \from X \to Z$ such that $x \mathrel{F} y \iff \pi(x) =
\pi(y)$ for all $x, y \in X$.

\Eigen-\Hajian-\Nadkarni negatively answered \Nadkarni's question
by providing a standard \Borel space $X$, a \Borel automorphism
$T \from X \to X$, and a smooth \Borel superequivalence relation
$F$ of the orbit equivalence relation $\orbitequivalencerelation{T}
{X}$ such that there is no weakly-wandering $T$-complete \Borel
set, but for each $F$-class $C$ there is a weakly-wandering
$(\restriction{T}{C})$-complete \Borel set (see \cite
{EHN93}). To see that the latter condition rules out the
existence of a $T$-invariant \Borel probability measure, observe that
the ergodic decomposition theorem ensures that the existence of a
$T$-invariant \Borel probability measure yields the existence of a
$T$-ergodic $T$-invariant \Borel probability measure (see \cite
[Theorem 3.2]{LM97} for the generalization to analytic
equivalence relations), the smoothness of $F$ implies that every
$T$-ergodic \Borel measure concentrates on some $F$-class $C$,
and the existence of a weakly-wandering $(\restriction
{T}{C})$-complete \Borel set rules out the existence of a
$(\restriction{T}{C})$-invariant \Borel probability measure.

However, \Eigen-\Hajian-\Nadkarni also noted that this leaves open
the question as to whether the inexistence of such an equivalence
relation yields a $T$-invariant \Borel probability measure.

The \definedterm{$T$-saturation} of a set $Y \subseteq X$ is given
by $\saturation{Y}{T} = \union[n \in \Z][\image{T^n}{Y}]$. When $X$
is a topological space, a homeomorphism $T \from X \to X$ is 
\definedterm{topologically transitive} if for all non-empty open sets
$U, V \subseteq X$ there exists $n \in \Z$ such that $\image{T^n}
{U} \intersection V \neq \emptyset$, and \definedterm{minimal} if
every orbit is dense. The \definedterm{odometer} is the isometry of
$\Cantorspace$ given by $\odometer(\constantsequence{1}{n}
\concatenation \sequence{0} \concatenation c) = \constantsequence
{0}{n} \concatenation \sequence{1} \concatenation c$ for all $c \in
\Cantorspace$ and $n \in \N$. Let $\muzero$ denote the \Borel
measure on $\Cantorspace$ given by $\muzero(\extensions{s}) = 1
/ 2^{\length{s}}$ for all $s \in \Cantortree$. It is easy to see that the
odometer is minimal, and therefore topologically transitive, and that
$\muzero$ is the unique $\sigma$-invariant \Borel probability
measure.

As noted in \cite[Example 3.13]{Mil04}, if $C \subseteq
\Cantorspace$ is a non-meager $\muzero$-null \Borel set and $X =
\saturation{C}{\odometer}$, then $\restriction{\sigma}{X}$ also yields
a negative answer to \Nadkarni's question. To see that there is no
weakly-wandering $(\restriction{\sigma}{X})$-complete \Borel set,
observe that if $B \subseteq \Cantorspace$ is a non-meager set with
the \Baire property, then there is a non-empty open set $U \subseteq
\Cantorspace$ in which $B$ is comeager (see, for example, \cite
[Proposition 8.26]{Kec95}), and the fact that $\sigma$ is a
homeomorphism ensures that the set $M = \saturation{U \setminus
B}{\sigma}$ is meager (see, for example, \cite[Exercise 8.45]
{Kec95}). Note also that the compactness of $\Cantorspace$ and
minimality of $\odometer$ yield a finite set $Z \subseteq \Z$ such
that $\Cantorspace = \union[n \in Z][\image{\odometer^n}{U}]$, in
which case $\Cantorspace = \union[n \in Z][\image{\odometer^n}
{\image{\odometer^k}{U}}]$ for all $k \in \Z$. In particular, it follows
that if $x \in X \setminus M$, then for all $k \in \Z$ there exists $n \in
Z$ such that $\odometer^{-n}(x) \in \image{\odometer^k}{U}$. Letting
$\calO$ denote the orbit of $x$ under $\sigma$, the pigeon-hole
principle therefore ensures that there is no set $S \subseteq \Z$ of
cardinality strictly greater than that of $Z$ for which $U \intersection
\calO$ is $S$-wandering, thus the same holds of $B \intersection
\calO$.

The transformation $\restriction{\odometer}{X}$ also yields a
negative answer to the subsequent question of
\Eigen-\Hajian-\Nadkarni, as the topological transitivity of
$\odometer$ ensures that every $\sigma$-invariant \Borel set is
comeager or meager, thus every smooth \Borel superequivalence
relation of $\orbitequivalencerelation{\odometer}{X}$ has a
comeager equivalence class.

Following \cite{Tse15}, we say that a set $Y \subseteq X$ is
\definedterm{locally weakly-wander\-ing} if its intersection with each
orbit of $T$ is weakly wandering. As noted by \Tserunyan, the
transformation $\restriction{\sigma}{X}$ also yields a negative
answer to the question as to whether the inexistence of a
locally-weakly-wandering $T$-complete \Borel set yields a
$T$-invariant \Borel probability measure. This is a direct
consequence of the above proof that there is no
weakly-wandering $(\restriction{\sigma}{X})$-complete \Borel set.

We say that a set $Y \subseteq X$ is \definedterm{very weakly
wandering} if there are arbitrarily large finite sets $S \subseteq \Z$
for which it is $S$-wandering. The above arguments also yield
negative answers to the analogous questions in which weak
wandering is replaced with very weak wandering.

Here we generalize these observations from \Borel actions of $\Z$
to \Borel actions of locally compact \Polish groups, and from \Borel
probability measures that are invariant with respect to a \Borel
action to \Borel probability measures that are invariant with respect
to a \Borel cocycle. We also show that not only are there examples
of \Borel actions yielding negative answers to the generalizations of
the questions considered by \Nadkarni, \Eigen-\Hajian-\Nadkarni,
and \Tserunyan, but that it was wholly unnecessary to search for
them in the first place, as they lie within every \Borel action that
could possibly contain them. Moreover, rather than just establishing
this for variants of weak wandering, we show that no recurrence
condition whatsoever yields the existence of an invariant \Borel
probability measure.

In \S\ref{recurrence}, we establish the basic properties of the
\definedterm{recurrence spectrum} of a \Borel action of a \Polish
group on a standard \Borel space, which codifies the suitably
robust forms of recurrence that it satisfies.

In \S\ref{maximal}, we show that locally-compact non-compact
\Polish groups have free \Borel actions on \Polish spaces with
maximal recurrence spectra.

In \S\ref{transience}, we show that the existence of
weakly-wandering and very-weakly-wandering suitably-complete
\Borel sets, as well as suitably-complete \Borel sets satisfying the
minimal non-trivial notion of transience corresponding to the failure
of the strongest notion of recurrence, in addition to the
\Eigen-\Hajian-\Nadkarni-style refinements thereof, can be
characterized in terms of the recurrence spectrum. Our arguments
also yield complexity bounds leading to implications between many
of these notions. For instance, it follows that if $X$ is a standard
\Borel space, $T \from X \to X$ is a \Borel automorphism, and there
is no smooth \Borel superequivalence relation $F$ of
$\orbitequivalencerelation{T}{X}$ with the property that there is a
weakly-wandering $(\restriction{T}{C})$-complete \Borel set for
every $F$-class $C$, then there is no locally-weakly-wandering
$T$-complete \Borel set.

In \S\ref{measures}, we generalize the generic compressibility
theorem of \Kechris-\Miller (see \cite[Theorem 13.1]{KM04})
to \Borel actions of locally compact \Polish groups on standard
\Borel spaces. We simultaneously replace comeagerness with a
stronger notion under which the recurrence spectrum is invariant,
thereby insuring that no condition on the latter yields an invariant
\Borel probability measure.

\section{The recurrence spectrum} \label{recurrence}
Suppose that $G \action X$ is a group action. We say that a set
$Y \subseteq X$ is \definedterm{complete} if $X = GY$, and
\definedterm{$\sigma$-complete} if there is a countable set $H
\subseteq G$ for which $X = HY$. The following observation
ensures that, under mild hypotheses, these notions coincide on
open sets.

\begin{proposition} \label{recurrence:stronglycomplete}
  Suppose that $G$ is a topological group, $H \subseteq G$ is
  dense, $X$ is a topological space, $G \action X$ is
  continuous-in-$G$, and $U \subseteq X$ is open. Then $GU =
  HU$.
\end{proposition}

\begin{propositionproof}
  Suppose that $g \in G$ and $x \in U$. As $G \action X$ is
  continuous-in-$G$, there is an open neighborhood $V \subseteq
  G$ of $\identity{G}$ for which $Vx \subseteq U$. Fix an open
  neighborhood $W \subseteq G$ of $g$ for which $\inverse{W} g
  \subseteq V$. As $H$ is dense, there exists $h \in H \intersection
  W$, and it only remains to observe that $\inverse{h} g \in V$, so
  $\inverse{h} g \cdot x \in U$, thus $g \cdot x \in hU$.
\end{propositionproof}

Let $\Ezero$ denote the equivalence relation on $\Cantorspace$
given by $c \mathrel{\Ezero} d \iff \exists n \in \N \forall m \ge n
\ c(m) = d(m)$. We will abuse language by saying that a subset
of $X$ is \definedterm{$\aleph_0$-universally \Baire} if its
pre-image under every \Borel function from a \Polish space to
$X$ has the \Baire property, and an
$\aleph_0$-universally-\Baire equivalence relation $E$ on $X$
is \definedterm{smooth} if there is no \Borel function $\pi \from
\Cantorspace \to X$ such that $c \mathrel{\Ezero} d \iff \pi(c)
\mathrel{E} \pi(d)$ for all $c, d \in \Cantorspace$. The
\Harrington-\Kechris-\Louveau generalization of the
\Glimm-\Effros dichotomy (see \cite[Theorem 1.1]
{HKL90}) ensures that this is compatible
with the usual notion of smoothness for \Borel equivalence
relations on standard \Borel spaces.

When $X$ is a topological space, a continuous-in-$X$ action $G
\action X$ is \definedterm{topologically transitive} if for all non-empty
open sets $U, V \subseteq X$ there exists $g \in G$ such that $gU
\intersection V \neq \emptyset$, and \definedterm{minimal} if every
orbit is dense.

\begin{proposition} \label{recurrence:universallyBaire}
  Suppose that $G$ is a group, $X$ is a \Polish space, $G \action
  X$ is continuous-in-$X$, $B \subseteq X$ is
  $\orbitequivalencerelation{G}{X}$-invariant, $C \subseteq X$ is an
  $\orbitequivalencerelation{G}{X}$-invariant \Gdelta set for which $G
  \action C$ is topologically transitive and in which $B$ is comeager,
  and $F$ is a smooth $\aleph_0$-universally-\Baire
  superequivalence relation of $\orbitequivalencerelation{G}{B}$.
  Then there is an $F$-class that is comeager in $C$.
\end{proposition}

\begin{propositionproof}
  Fix a dense \Gdelta set $C' \subseteq C$ contained in $B$, and
  note that $F$ has the \Baire property in $C' \times C'$, thus
  in $C \times C$. The straightforward generalization of the
  \Becker-\Kechris criterion for continuously embedding $\Ezero$
  from orbit equivalence relations induced by groups of
  homeomorphisms (see \cite[Theorem 3.4.5]{BK96}) to
  superequivalence relations of such orbit equivalence relations (see,
  for example, \cite[Theorem 2.1]{KMS14}) ensures that
  the union of $F$ and $(C \setminus B) \times (C \setminus B)$ is
  non-meager in $C \times C$, so the \Kuratowski-\Ulam theorem
  (see, for example, \cite[Theorem 8.41]{Kec95}) yields an
  $F$-class that is non-meager and has the \Baire property in $C$,
  thus comeager in $C$ by topological transitivity.
\end{propositionproof}

For each set $Y \subseteq X$, define $\returntimes{Y} = \set{g \in
G}[Y \intersection gY \neq \emptyset]$. The following fact is the
obvious generalization of \Pettis's Lemma (see, for example, \cite
[Theorem 9.9]{Kec95}) to group actions.

\begin{proposition} \label{recurrence:Pettis:global}
  Suppose that $G$ is a group, $X$ is a \Baire space, $G \action
  X$ is continuous-in-$X$, $U \subseteq X$ is non-empty and
  open, and $B \subseteq U$ is comeager. Then $\returntimes
  {U} \subseteq \returntimes{B}$.
\end{proposition}

\begin{propositionproof}
  Note that for all $g \in G$, the fact that $G \action X$ is
  continuous-in-$X$ ensures that $gU$ is open and $gB$ is
  comeager in $gU$. In particular, it follows that if $U \intersection
  gU$ is non-empty, then $B \intersection gB$ is comeager in $U
  \intersection gU$, in which case the fact that $X$ is a \Baire space
  ensures that $B \intersection gB$ is also non-empty, thus
  $\returntimes{U} \subseteq \returntimes{B}$.
\end{propositionproof}

Given an upward-closed family $\calF \subseteq \powerset{G}$ and
any family $\Gamma \subseteq \powerset{X}$, we say that an
action $G \action X$ is \definedterm{$\calF$-recurrent on
$\Gamma$} if $\returntimes{B} \in \calF$ for all $B \in \Gamma$.
This generalizes the usual notion of recurrence in topological
dynamics, where one says that a continuous-in-$X$ action of a
group on a topological space $X$ is \definedterm{$\calF$-recurrent} if
it is $\calF$-recurrent on non-empty open sets. We next note that,
under mild hypotheses, this latter notion propagates to recurrence
on non-meager sets with the \Baire property.

\begin{proposition} \label{recurrence:openBaire:global}
  Suppose that $G$ is a group, $X$ is a \Baire space, $\calF
  \subseteq \powerset{G}$ is upward closed, and $G \action X$ is
  continuous-in-$X$ and $\calF$-recurrent. Then $G \action X$ is
  $\calF$-recurrent on non-meager sets with the \Baire property.
\end{proposition}

\begin{propositionproof}
  If $B \subseteq X$ is a non-meager set with the \Baire property,
  then there is a non-empty open set $U \subseteq X$ in which $B$
  is comeager, and Proposition \ref{recurrence:Pettis:global}
  ensures that $\returntimes{U} \subseteq \returntimes{B}$, thus
  $\returntimes{B} \in \calF$.
\end{propositionproof}

When $X$ is a \Polish space, the \definedterm{decomposition
into minimal components} of a continuous action $G \action X$
is the equivalence relation on $X$ given by $x \mathrel
{\minimaldecomposition{G}{X}} y \iff \closure{\orbit{x}{G}} = \closure
{\orbit{y}{G}}$. It is easy to see that $\minimaldecomposition{G}{X}$
is \Gdelta and smooth (although the latter is also a direct consequence
of the former and \cite[Corollary 1.2]{HKL90}), and
that for each $\minimaldecomposition{G}{X}$-class $C$, the action $G
\action C$ is minimal.

When $G$ is \Borel, the \definedterm{recurrence spectrum} of a
\Borel action $G \action X$ is the collection of all upward-closed
families $\calF \subseteq \powerset{G}$ such that every smooth
\Borel superequivalence relation $F$ of $\orbitequivalencerelation
{G}{X}$ has an equivalence class $C$ for which $G \action C$ is
$\calF$-recurrent on $\sigma$-complete \Borel sets.  As a theorem
of \Becker-\Kechris ensures that every \Borel action of a \Polish
group on a standard \Borel space is \Borel isomorphic to a
continuous action on a \Polish space (see \cite[Theorem 5.2.1]
{BK96}), the following observation ensures that, under
mild hypotheses, the notion of recurrence spectrum is robust, in
the sense that it does not depend on the particular underlying notion
of definability, and in the sense that it is invariant under passage to
sufficiently large $\orbitequivalencerelation{G}{X}$-invariant subsets.

\begin{proposition} \label{recurrence:characterization:global}
  Suppose that $G$ is a separable group, $X$ is a
  \Polish space, $G \action X$ is continuous, $\calF \subseteq
  \powerset{G}$ is upward closed, and $B \subseteq X$ is
  $\orbitequivalencerelation{G}{X}$-invariant and comeager in every
  $\minimaldecomposition{G}{X}$-class. Then the following are
  equivalent:
  \begin{enumerate}
    \item Every smooth $\aleph_0$-universally-\Baire
    superequivalence relation $F$ of $\orbitequivalencerelation{G}
      {B}$ has a class $C$ for which $G \action C$ is
      $\calF$-recurrent on $\sigma$-complete
      $\aleph_0$-universally-\Baire sets.
    \item There is an $\minimaldecomposition{G}{X}$-class $C$ for
      which $G \action C$ is $\calF$-recurrent.
  \end{enumerate}
\end{proposition}

\begin{propositionproof}
  To see $(1) \implies (2)$, fix an $\minimaldecomposition{G}
  {X}$-class $C$ for which $G \action B \intersection C$ is
  $\calF$-recurrent on $\sigma$-complete open sets. To see that $G
  \action C$ is $\calF$-recurrent, suppose that $U \subseteq C$ is a
  non-empty open set, and note that the minimality of $G \action C$
  ensures that $U$ is complete, and therefore $\sigma$-complete by
  Proposition \ref{recurrence:stronglycomplete}, thus $\returntimes{U}
  \in \calF$.
  
  To see $(2) \implies (1)$, fix an $\minimaldecomposition{G}
  {X}$-class $C$ for which $G \action C$ is $\calF$-recurrent, and
  suppose that $F$ is a smooth $\aleph_0$-universally-\Baire
  superequivalence relation of $\orbitequivalencerelation{G}{B}$.
  Proposition \ref{recurrence:universallyBaire} then yields an
  $F$-class $D$ that is comeager in $C$. To see that $G \action D$
  is $\calF$-recurrent on $\sigma$-complete
  $\aleph_0$-universally-\Baire sets, suppose that $A \subseteq D$
  is such a set, and note that $\sigma$-completeness ensures that
  $A$ is non-meager in $C$. Fix a dense \Gdelta set $C' \subseteq
  C$ contained in $D$, and note that $A \intersection C'$ has the
  \Baire property in $C'$, so $A$ has the \Baire property in $C$,
  thus Proposition \ref{recurrence:openBaire:global} ensures that
  $\returntimes{A} \in \calF$.
\end{propositionproof}

Let $\dual{\Gamma}$ denote the family of sets whose complements
are in $\Gamma$, let $\Dtwo{\Gamma}$ denote the family of
differences of sets in $\Gamma$, and let $\unions{\Gamma}$
denote the family of countable unions of sets in $\Gamma$.

The \definedterm{horizontal sections} of a set $R \subseteq X \times
Y$ are the sets of the form $\horizontalsection{R}{y} = \set{x \in X}
[x \mathrel{R} y]$ for $y \in Y$, whereas the \definedterm{vertical
sections} of a set $R \subseteq X \times Y$ are the sets of the form
$\verticalsection{R}{x} = \set{y \in Y}[x \mathrel{R} y]$ for $x \in
X$. We say that $\calF$ is \definedterm{$\Gamma$-on-open} if
$\set{x \in X}[\horizontalsection{U}{x} \in \calF] \in \Gamma$ for all
open sets $U \subseteq G \times X$.

Given a superequivalence relation $E$ of $\orbitequivalencerelation
{G}{X}$, we say that an action $G \action X$ is \definedterm
{$E$-locally $\calF$-recurrent on $\Gamma$} if for all $B \in
\Gamma$, there is an $E$-class $C$ such that $\returntimes{B
\intersection C} \in \calF$. We next note that, under mild hypotheses,
the recurrence spectrum can also be characterized in terms of local
recurrence of $G \action X$ itself.

\begin{proposition} \label{recurrence:complete}
  Suppose that $G$ is a separable group, $X$ is a \Polish space, $G
  \action X$ is continuous, $\Gamma \subseteq \powerset{X}$ is a
  family of $\aleph_0$-universally-\Baire sets containing the open
  sets and closed under finite intersections and finite unions, and
  $\calF \subseteq \powerset{G}$ is upward closed and $\dual
  {\Gamma}$-on-open. Then the following are equivalent:
  \begin{enumerate}
    \item There is an $\minimaldecomposition{G}{X}$-class $C$ for
      which $G \action C$ is $\calF$-recurrent.
    \item The action $G \action X$ is $\minimaldecomposition{G}
      {X}$-locally $\calF$-recurrent on $\sigma$-com\-plete
      $\unions{(\Dtwo{\Gamma})}$ sets.
  \end{enumerate}
\end{proposition}

\begin{propositionproof}
  To see $\neg (2) \implies \neg (1)$, observe that if $B \subseteq X$
  is a $\sigma$-complete $\unions{(\Dtwo{\Gamma})}$ set such
  that $\returntimes{B \intersection C} \notin \calF$ for every
  $\minimaldecomposition{G}{X}$-class $C$, then it is non-meager
  and has the \Baire property in every such class, so there is no
  $\minimaldecomposition{G}{X}$-class $C$ for which $G \action C$
  is $\calF$-recurrent by Proposition \ref{recurrence:openBaire:global}.

  To see $\neg (1) \implies \neg (2)$, fix a basis $\sequence{U_n}[n
  \in \N]$ for $X$. For all $n \in \N$, define $V_n = \set{\pair{g}{x} \in
  G \times X}[U_n \intersection gU_n \intersection \equivalenceclass
  {x}{\minimaldecomposition{G}{X}} \neq \emptyset]$. Observe
  that if $\pair{g}{x} \in V_n$, then the minimality of $G \action
  \equivalenceclass{x}{\minimaldecomposition{G}{X}}$ yields $h \in
  G$ for which $h \cdot x \in U_n \intersection gU_n$, so the
  continuity of $G \action X$ yields open neighborhoods $U_g
  \subseteq G$ of $g$ and $U_x \subseteq X$ of $x$ such that
  $hU_x \union \inverse{U_g}hU_x \subseteq U_n$, thus $U_g
  \times U_x \subseteq V_n$, hence $V_n$ is open. It follows that
  the $\minimaldecomposition{G}{X}$-invariant sets $A_n = \set{x
  \in G U_n}[\horizontalsection{V_n}{x} \notin \calF]$ are in
  $\Gamma$, so the sets $B_n = A_n \setminus \union[m < n]
  [A_m]$ are in $\Dtwo{\Gamma}$, thus the set $B = \union[n \in
  \N][B_n \intersection U_n]$ is in $\unions{(\Dtwo{\Gamma})}$.
  But if there is no $\minimaldecomposition{G}{X}$-class $C$ for
  which $G \action C$ is $\calF$-recurrent, then $B$ is complete,
  and therefore $\sigma$-complete by Proposition \ref
  {recurrence:stronglycomplete}, thus $G \action X$ is not
  $\minimaldecomposition{G}{X}$-locally $\calF$-recurrent on
  $\sigma$-complete $\unions{(\Dtwo{\Gamma})}$ sets.
\end{propositionproof}

We next show that, under an additional mild hypothesis on $G$,
Propositions \ref{recurrence:Pettis:global}, \ref
{recurrence:openBaire:global}, and \ref
{recurrence:characterization:global} can be strengthened so as to
show that the recurrence spectrum is also robust in the sense that
it does not depend on whether the underlying notion of recurrence
is local.

The following fact is a somewhat more intricate generalization of
the special case of \Pettis's Lemma for second-countable groups
to continuous actions of such groups.

\begin{proposition} \label{recurrence:Pettis:local}
  Suppose that $G$ is a second-countable \Baire group, $X$ is a
  second-countable \Baire space, $G \action X$ is continuous, $U
  \subseteq X$ is non-empty and open, and $B \subseteq U$ is
  comeager. Then $\returntimes{U \intersection \orbit{x}{G}}
  \subseteq \returntimes{B \intersection \orbit{x}{G}}$ for comeagerly
  many $x \in X$.
\end{proposition}

\begin{propositionproof}
  We write $\forcomeagerlymany x \in X \ \phi(x)$ to indicate that
  $\set{x \in X}[\phi(x)]$ is comeager. As the fact that $G \action X$
  is continuous-in-$X$ ensures that it is open, it follows that $\set
  {\pair{g}{x} \in G \times X}[g \cdot x \notin U \setminus B]$ is
  comeager, so the set $C = \set{x \in X}[\forcomeagerlymany g \in
  G \ g \cdot x \notin U \setminus B]$ is comeager by the
  \Kuratowski-\Ulam theorem. Observe that if $h \in G$ and $x \in X$,
  then $\set{g \in G}[g \cdot x \notin h (U \setminus B)] = h \set{g \in
  G}[g \cdot x \notin U \setminus B]$, so the fact that $G \action X$ is
  continuous-in-$X$ also ensures that if $x \in C$, then
  $\forcomeagerlymany g \in G \ g \cdot x \notin (U \setminus B)
  \union h(U \setminus B)$, in which case the fact that $(U
  \intersection hU) \setminus (B \intersection hB) \subseteq (U
  \setminus B) \union h(U \setminus B)$ therefore implies that
  $\forcomeagerlymany g \in G \ g \cdot x \notin (U \intersection hU)
  \setminus (B \intersection hB)$. Note now that for all $h \in G$,
  the fact that $G \action X$ is continuous-in-$X$ ensures that $U
  \intersection hU$ is open, so the fact that $G \action X$ is
  continuous-in-$G$ implies that if $x \in X$ and $U \intersection
  hU \intersection \orbit{x}{G}$ is non-empty, then there are
  non-meagerly many $g \in G$ for which $g \cdot x \in U
  \intersection hU$. In particular, it follows that if $x \in C$ and $U
  \intersection hU \intersection \orbit{x}{G}$ is non-empty, then so
  too is $B \intersection hB \intersection \orbit{x}{G}$, hence
  $\returntimes{U \intersection \orbit{x}{G}} \subseteq \returntimes
  {B \intersection \orbit{x}{G}}$ for all $x \in C$.
\end{propositionproof}

We next note that, under mild hypotheses, $\calF$-recurrence of
topologically transitive actions not only propagates to
$\calF$-recurrence on non-meager sets with the \Baire property,
but to its $\orbitequivalencerelation{G}{X}$-local strengthening.

\begin{proposition} \label{recurrence:openBaire:local}
  Suppose that $G$ is a second-countable \Baire group, $X$ is a
  second-countable \Baire space, $\calF \subseteq \powerset{G}$ is
  upward closed, and $G \action X$ is continuous, $\calF$-recurrent,
  and topologically transitive. Then $G \action X$ is
  $\orbitequivalencerelation{G}{X}$-locally $\calF$-recurrent on
  non-meager sets with the \Baire property.
\end{proposition}

\begin{propositionproof}
  Suppose that $B \subseteq X$ is a non-meager set with the \Baire
  property, and fix a non-empty open set $U \subseteq X$ in
  which $B$ is comeager. The topological transitivity of $G \action X$
  ensures that the set $C = \set{x \in X}[\orbit{x}{G} \text{ is dense}]$
  is comeager, and Proposition \ref{recurrence:Pettis:local} implies
  that the set $D = \set{x \in X}[\returntimes{U \intersection \orbit{x}{G}}
  \subseteq \returntimes{B \intersection \orbit{x}{G}}]$ is comeager.
  So it only remains to observe that if $x \in C \intersection D$, then
  $\returntimes{U} \subseteq \returntimes{U \intersection \orbit{x}{G}}
  \subseteq \returntimes{B \intersection \orbit{x}{G}}$, thus
  $\returntimes{B \intersection \orbit{x}{G}} \in \calF$.
\end{propositionproof}

We can now establish the promised robustness result.

\begin{proposition} \label{recurrence:characterization:local}
  Suppose that $G$ is a second-countable \Baire group, $X$ is a
  \Polish space, $G \action X$ is continuous, $\calF \subseteq
  \powerset{G}$ is upward closed, and $B \subseteq X$ is
  $\orbitequivalencerelation{G}{X}$-invariant and comeager in
  every $\minimaldecomposition{G}{X}$-class. Then the following
  are equivalent:
  \begin{enumerate}
    \item Every smooth $\aleph_0$-universally-\Baire superequivalence
      relation $F$ of $\orbitequivalencerelation{G}{B}$ has a class $C$
      for which $G \action C$ is $\orbitequivalencerelation{G}
      {C}$-locally $\calF$-recurrent on $\sigma$-complete
      $\aleph_0$-universally-\Baire sets.
    \item There is an $\minimaldecomposition{G}{X}$-class $C$ for
      which $G \action C$ is $\calF$-recurrent.
  \end{enumerate}
\end{proposition}

\begin{propositionproof}
  Exactly as in the proof of Proposition \ref
  {recurrence:characterization:global}, but replacing the use of
  Proposition \ref{recurrence:openBaire:global}
  with that of Proposition \ref{recurrence:openBaire:local}.
\end{propositionproof}

\section{The strongest notion of recurrence} \label{maximal}

Given a set $S \subseteq G$, we say that a set $Y \subseteq X$ is
\definedterm{$S$-transient} if $Y \intersection SY = \emptyset$.

\begin{proposition} \label{maximal:cocompact:transient}
  Suppose that $G$ is a topological group, $X$ is a \Hausdorff
  space, $G \action X$ is continuous, $K \subseteq G$ is compact,
  and $x \in X$ is not fixed by any element of $K$. Then there is a
  $K$-transient open neighborhood of $x$.
\end{proposition}

\begin{propositionproof}
  For each $g \in K$, the fact that $X$ is \Hausdorff yields disjoint
  open neighborhoods $V_g \subseteq X$ of $x$ and $W_g
  \subseteq X$ of $g \cdot x$, and the continuity of $G \action X$
  yields open neighborhoods $U_g \subseteq G$ of $g$ and $V_g'
  \subseteq V_g$ of $x$ for which $U_g V_g' \subseteq W_g$,
  thus $U_g V_g' \intersection V_g' = \emptyset$. The compactness
  of $K$ then yields a finite set $F \subseteq K$ for which $K
  \subseteq \union[g \in F][U_g]$, in which case $\intersection[g \in
  F][V_g']$ is a $K$-transient open neighborhood of $x$.
\end{propositionproof}

It follows that upward-closed families whose corresponding notions
of recurrence are realizable by suitable free actions necessarily
contain all co-compact neighborhoods of the identity.

\begin{proposition}
  Suppose that $G$ is a topological group and $\calF \subseteq
  \powerset{G}$ is an upward closed family for which there is an
  $\calF$-recurrent continuous free action of $G$ on a \Hausdorff
  space. Then $\calF$ contains every co-compact neighborhood
  of $\identity{G}$.
\end{proposition}

\begin{propositionproof}
  This is a direct consequence of Proposition \ref
  {maximal:cocompact:transient}.
\end{propositionproof}

Given sets $Y, Z \subseteq X$, define $\hittingtimes{Y}{Z} =
\set{g \in G}[gY \intersection Z \neq \emptyset]$. When $X$ is a
topological space, we say that a continuous-in-$X$ action $G
\action X$ is \definedterm{topologically mixing} if $\hittingtimes
{U}{V}$ is co-compact for all non-empty open sets $U, V \subseteq
X$.

\begin{proposition} \label{examples:cocompact:topological}
  Suppose that $G$ is a topological group, $X$ is a topological
  space, $\calF \subseteq \powerset{G}$ is the family of
  co-pre-compact subsets of $G$ containing $\identity{G}$, and
  $G \action X$ is continuous-in-$X$, $\calF$-recurrent, and
  topologically transitive. Then $G \action X$ is topologically
  mixing.
\end{proposition}

\begin{propositionproof}
  Given non-empty open sets $U, V \subseteq X$, the topological
  transitivity of $G \action X$ yields $g \in G$ for which $gU
  \intersection V \neq \emptyset$, so the fact that $G
  \action X$ is $\calF$-recurrent ensures that $\returntimes{gU
  \intersection V}$ is co-compact. As $h \in \returntimes
  {gU \intersection V} \implies hgU \intersection V \neq \emptyset
  \implies hg \in \hittingtimes{U}{V}$, it follows that $\returntimes{gU
  \intersection V} g \subseteq \hittingtimes{U}{V}$, so $\hittingtimes
  {U}{V}$ is co-compact.
\end{propositionproof}

It follows that the existence of a suitable free \Borel action of $G$
whose recurrence spectrum contains the family $\calF$ of
co-pre-compact subsets of $G$ containing $\identity{G}$ is
equivalent to the existence of a suitable continuous
topologically-mixing free action of $G$.

\begin{proposition}
  Suppose that $G$ is a separable group, $X$ is a \Polish
  space, $G \action X$ is continuous, and $\calF \subseteq
  \powerset{G}$ is the family of co-pre-compact subsets of $G$
  containing $\identity{G}$. Then $\calF$ is in the recurrence
  spectrum of $G \action X$ if and only if there is an equivalence
  class $C$ of $\minimaldecomposition{G}{X}$ for which $G
  \action C$ is topologically mixing.
\end{proposition}

\begin{propositionproof}
  This is a direct consequence of Propositions \ref
  {recurrence:characterization:global} and \ref
  {examples:cocompact:topological}.
\end{propositionproof}

To our suprise, we were unable to find a proof in the literature of
the fact that locally-compact non-compact \Polish groups have
free topologically-mixing continuous actions on \Polish spaces. In a
pair of private emails, \Glasner-\Weiss suggested that the
strengthening in which the underlying space is compact should be a
consequence of generalizations of the results of \cite{Wei12} to
locally compact groups, and that a substantially simpler construction
should yield the aforementioned result. However, we give an
elementary proof by checking that the action of $G$ by left
multiplication on the space $\closedsets{G}$ of closed subsets of
$G$ is topologically mixing, where $\closedsets{G}$ is equipped
with the \definedterm{\Fell topology} generated by the sets $V_K =
\set{F \in \closedsets{G}}[F \intersection K = \emptyset]$ and
$W_U = \set{F \in \closedsets{G}}[F \intersection U \ne \emptyset]$,
where $K \subseteq G$ compact and $U \subseteq G$ open. It is
well-known that $\closedsets{G}$ is a compact \Polish space (see,
for example, \cite[Exercise 12.7]{Kec95}).

\begin{proposition}
  Suppose that $G$ is a locally-compact non-compact \Polish group.
  Then there is a \Polish space $X$ for which there is a free
  topologically-mixing continuous action $G \action X$.
\end{proposition}

\begin{propositionproof}
  While it is well-known that $G \action \closedsets{G}$ is
  continuous, we will provide a proof for the reader's convenience.
  Towards this end, it is sufficient to show that if $g \in G$, $F \in
  \closedsets{G}$, and $U_{gF} \subseteq \closedsets{G}$ is an
  open neighborhood of $gF$, then there are open neighborhoods
  $U_g \subseteq G$ of $g$ and $U_F \subseteq \closedsets{G}$ of
  $F$ for which $U_g U_F \subseteq U_{gF}$. Clearly we can
  assume that $U_{gF} = V_K$ for some compact $K \subseteq G$,
  or $U_{gF} = W_U$ for some open set $U \subseteq G$. In the
  former case, it follows that $F \intersection \inverse{g}K =
  \emptyset$, so the local compactness of $G$ ensures that for all
  $h \in K$ there are a pre-compact open neighborhood $U_{g,h}
  \subseteq G$ of $g$ and an open neighborhood $V_{g,h} \subseteq
  G$ of $h$ such that $F \intersection \inverse{\closure{U_{g,h}}}
  V_{g,h} = \emptyset$, and the compactness of $K$ yields a finite
  set $L \subseteq K$ such that $K \subseteq \union[h \in L][V_{g,h}]$,
  in which case the sets $U_g = \intersection[h \in L][U_{g,h}]$ and
  $U_F = V_{\inverse{\closure{U_g}} K}$ are as desired. In the latter
  case, there exists $h \in F$ for which $gh \in U$, so there are open
  neighborhoods $U_g \subseteq G$ of $g$ and $U_h \subseteq G$
  of $h$ such that $U_g U_h \subseteq U$, thus the sets $U_g$ and
  $U_F = W_{U_h}$ are as desired.
  
  To see that $G \action \closedsets{G}$ is topologically mixing, it is
  sufficient to show that if $U = V_K \intersection \intersection[i < m]
  [W_{U_i}]$ and $U' = V_{K'} \intersection \intersection[j < n]
  [W_{U_j'}]$ are non-empty, where $K, K' \subseteq G$ are
  compact and $U_i \subseteq \setcomplement{K}$ and $U_j'
  \subseteq \setcomplement{K'}$ are open for all $i < m$ and $j <
  n$, then $\set{g \in G}[g U \intersection U' = \emptyset]$ is
  compact. Towards this end, note that for all $i < m$ and $j < n$,
  the sets
  \begin{equation*}
    \textstyle
    L_i
      = \set{g \in G}[g U_i \subseteq K']
      = \intersection[h \in U_i][{\set{g \in G}
        [g h \in K']}]
      = \intersection[h \in U_i][K' \inverse{h}]
  \end{equation*}
  \centerline{and}
  \begin{equation*}
    \textstyle
    L_j'
      = \set{g \in G}[U_j' \subseteq g K]
      = \intersection[h \in U_j'][{\set{g \in G}
        [h \in g K]}]
      = \intersection[h \in U_j'][h \inverse{K}]
  \end{equation*}
  are compact, set $L = \union[i < m][L_i] \union \union[j < n][L_j']$,
  and observe that if $g \in \setcomplement{L}$, then there exist
  $g_i \in gU_i \setminus K'$ for all $i < m$ and $g_j' \in U_j'
  \setminus gK$ for all $j < n$. But then the set $F = \set{g_i}[i
  < n] \union \set{g_j'}[j < n]$ is in $gU \intersection U'$.
  
  The \definedterm{free part} of the action $G \action \closedsets{G}$
  is the set $X$ of $F \in \closedsets{G}$ that are not fixed by any
  non-identity element of $G$. The local-compactness and
  separability of $G$ ensure that $X$ is the intersection of
  countably-many sets of the form $X_K = \set{F \in \closedsets{G}}
  [\forall g \in K \ gF \neq F]$, where $K \subseteq G \setminus
  \set{\identity{G}}$ is compact. As Proposition \ref
  {maximal:cocompact:transient} ensures that each $X_K$ is open,
  it follows that $X$ is \Gdelta, and therefore \Polish. To see that $G
  \action X$ is the desired action, it only remains to establish that
  $X$ is comeager. And for this, it is sufficient to show that if $K
  \subseteq G \setminus \set{\identity{G}}$ is compact, then $X_K$
  is dense. Towards this end, suppose that $U = V_L \intersection
  \intersection[i < n][W_{U_i}]$ is non-empty, where $L \subseteq G$
  is compact and $U_i \subseteq \setcomplement{L}$ is open for all
  $i < n$, and fix $g_i \in U_i$ for all $i < n$. As $G$ is locally compact,
  by passing to open neighborhoods of $g_i$ contained in $U_i$, we
  can assume that each of the sets $U_i$ is pre-compact. As $G$ is
  not compact, there exists $g \in \setcomplement{(L \union \union[i
  < n][\inverse{K} U_i])}$. Then the set $F = \set{g} \union \set{g_i}[i
  < n]$ is in $U$, and the fact that $F \intersection Kg = \emptyset$
  ensures that $F \in X_K$.
\end{propositionproof}

\section{Transience} \label{transience}

Following the usual convention, let $\Sigmaclass[0][1]$ denote
the pointclass of open sets, and recursively define $\Piclass[0][n]$
to be the pointclass of complements of $\Sigmaclass[0][n]$ sets,
and $\Sigmaclass[0][n+1]$ to be the pointclass of countable unions
of $\Piclass[0][n]$ sets. It is well known that $\unions{(\Dtwo
{\Sigmaclass[0][n]})}$, $\unions{(\Dtwo{\Piclass[0][n]})}$, and
$\Sigmaclass[0][n+1]$ coincide on metric spaces (see,
for example, \cite[\S11.B]{Kec95}).

Given a set $S \subseteq G$, we use $\principleF{S}$ to denote
the family of sets $T \subseteq G$ for which $S \intersection T
\neq \emptyset$. 

\begin{proposition} \label{transience:Stransient:global}
  Suppose that $G$ is a separable group, $X$ is a \Polish
  space, $G \action X$ is continuous, and $S \subseteq G$. Then
  the following are equivalent:
  \begin{enumerate}
    \item The family $\principleF{S}$ is not in the recurrence spectrum
      of $G \action X$.
    \item There is a smooth $\aleph_0$-universally \Baire
      superequivalence $F$ of $\orbitequivalencerelation{G}{X}$
      for which each action $G \action \equivalenceclass{x}{F}$
      has an $S$-transient $\sigma$-complete
      $\aleph_0$-universally-\Baire set.
    \item The action $G \action X$ has an $S$-transient
      $\sigma$-complete $\Sigmaclass[0][2]$ set.
  \end{enumerate}
\end{proposition}

\begin{propositionproof}
  As a set $Y \subseteq X$ is $S$-transient if and only if
  $\returntimes{Y} \notin \principleF{S}$, Proposition \ref
  {recurrence:characterization:global} yields $(1) \iff (2)$. Note that if
  $U \subseteq G \times X$ and $x \in X$, then $\horizontalsection
  {U}{x} \in \principleF{S} \iff x \in \union[g \in S][\verticalsection{U}
  {g}]$, so $\principleF{S}$ is $\Sigmaclass[0][1]$-on-open. As a
  set $Y \subseteq X$ is $S$-transient if and only if $\returntimes{C
  \intersection Y} \notin \principleF{S}$ for all equivalence classes $C$
  of $\minimaldecomposition{G}{X}$, Proposition \ref
  {recurrence:complete} yields $(1) \iff (3)$.
\end{propositionproof}

Given a set $\calS \subseteq \powerset{G}$, we say that a set $Y
\subseteq X$ is \definedterm{$\calS$-transient} if there is a set $S
\in \calS$ for which $Y$ is $S$-transient.

\begin{proposition} \label{transience:calStransient:global}
  Suppose that $G$ is a separable group, $X$ is a \Polish
  space, $G \action X$ is continuous, and $\calS \subseteq
  \powerset{G}$. Then the following are equivalent:
  \begin{enumerate}
    \item There exists $S \in \calS$ for which the family $\principleF
      {S}$ is not in the recurrence spectrum of $G \action X$.
    \item There exist $S \in \calS$ and a smooth
      $\aleph_0$-universally \Baire superequivalence $F$ of
      $\orbitequivalencerelation{G}{X}$ for which each action $G
      \action \equivalenceclass{x}{F}$ has an $S$-transient $\sigma$-complete
      $\aleph_0$-universally-\Baire set.
    \item The action $G \action X$ has an $\calS$-transient $\sigma$-complete
      $\Sigmaclass[0][2]$ set.
  \end{enumerate}
\end{proposition}

\begin{propositionproof}
  This is a direct consequence of Proposition \ref
  {transience:Stransient:global}.
\end{propositionproof}

Given a set $S \subseteq G$, we say that a set $Y \subseteq X$ is
\definedterm{$S$-wandering} if $gY \intersection hY = \emptyset$
for all distinct $g, h \in S$, and \definedterm{weakly wandering} if
there exists an infinite set $S \subseteq G$ for which $Y$ is
$S$-wandering.

\begin{proposition} \label{transience:weakwandering:global}
  Suppose that $G$ is a separable group, $X$ is a \Polish space,
  $G \action X$ is continuous, and $\calS$ is the family of sets of
  the form $\inverse{S} S \setminus \set{\identity{G}}$, where $S
  \subseteq G$ is infinite. Then the following are equivalent:
  \begin{enumerate}
    \item There exists $S \in \calS$ for which the family $\principleF
      {S}$ is not in the recurrence spectrum of $G \action X$.
    \item There exist an infinite set $S \subseteq G$ and a smooth
      $\aleph_0$-universally \Baire superequivalence $F$ of
      $\orbitequivalencerelation{G}{X}$ for which each action $G
      \action \equivalenceclass{x}{F}$ has an $S$-wandering
      $\sigma$-complete $\aleph_0$-universally-\Baire set.
    \item The action $G \action X$ has a weakly-wandering
      $\sigma$-complete $\Sigmaclass[0][2]$ set.  
  \end{enumerate}
\end{proposition}

\begin{propositionproof}
  Observe that if $S \subseteq G$, then a set $Y \subseteq X$ is
  $S$-wandering if and only if it is $(\inverse{S} S \setminus
  \set{\identity{G}})$-transient, and appeal to Proposition \ref
  {transience:calStransient:global}.
\end{propositionproof}

A subset of a standard \Borel space is \definedterm{analytic} if it is
the image of a standard \Borel space under a \Borel function, and
\definedterm{co-analytic} if its complement is analytic. A result of
\Lusin-\Sierpinski (see, for example, \cite[Theorem 21.6]{Kec95})
ensures that analytic subsets of standard \Borel spaces are
$\aleph_0$-universally \Baire. Let $\Sigmaclass[1][1]$ denote the
pointclass of analytic sets, and let $\Piclass[1][1]$ denote the
pointclass of co-analytic sets.

Given a superequivalence relation $E$ of $\orbitequivalencerelation
{G}{X}$, we say that a set $Y \subseteq X$ is \definedterm
{$E$-locally weakly-wandering} if its intersection with each
$E$-class is weakly wandering. Given a set $\calS \subseteq
\powerset{G}$, define $\principleF{\calS} = \intersection[S \in \calS]
[\principleF{S}]$.

\begin{proposition} \label{transience:weakwandering:local}
  Suppose that $G$ is a \Polish group, $X$ is a \Polish space, $G
  \action X$ is continuous, and $\calS$ is the family of sets of the
  form $\inverse{S} S \setminus \set{\identity{G}}$, where $S
  \subseteq G$ is infinite. Then the following are equivalent:
  \begin{enumerate}
    \item The family $\principleF{\calS}$ is not in the
      recurrence spectrum of $G \action X$.
    \item There is a smooth $\aleph_0$-universally \Baire
      superequivalence $F$ of $\orbitequivalencerelation{G}{X}$ for
      which each action $G \action \equivalenceclass{x}{F}$ has an
      $\orbitequivalencerelation{G}{X}$-locally-weakly-wandering
      $\sigma$-complete $\aleph_0$-universally-\Baire set.
    \item The action $G \action X$ has an $\minimaldecomposition
      {G}{X}$-locally-weakly-wandering $\sigma$-complete $\unions
      {(\Dtwo{\Sigmaclass[1][1]})}$ set.
  \end{enumerate}
\end{proposition}

\begin{propositionproof}
  As a set $Y \subseteq X$ is $\orbitequivalencerelation{G}
  {X}$-locally weakly-wandering if and only if $\returntimes{C
  \intersection Y} \notin \principleF{\calS}$ for all equivalence classes
  $C$ of $\orbitequivalencerelation{G}{X}$, Proposition \ref
  {recurrence:characterization:local} yields $(1) \iff (2)$.
  Note that if $U \subseteq G \times X$ and $x \in X$,
  then $\horizontalsection{U}{x} \in \calF_{\calS} \iff \forall
  \sequence{g_i}[i \in \N] \in \functions{\N}{G} \exists i \neq j
  \ ( g_i = g_j \mathor \inverse{g_i} g_j \in \horizontalsection
  {U}{x})$, so $\principleF{\calS}$ is $\Piclass[1][1]$-on-open. As a
  set $Y \subseteq X$ is $\minimaldecomposition{G}{X}$-locally
  weakly-wandering if and only if $\returntimes{C \intersection Y}
  \notin \principleF{\calS}$ for all equivalence classes $C$ of
  $\minimaldecomposition{G}{X}$, Proposition \ref
  {recurrence:complete} yields $(1) \iff (3)$.
\end{propositionproof}

We say that a set $Y \subseteq X$ is \definedterm{very weakly
wandering} if there are arbitrarily large finite sets $S \subseteq G$
for which $Y$ is $S$-wandering.

\begin{proposition} \label{transience:veryweakwandering:global}
  Suppose that $G$ is a separable group, $X$ is a \Polish space,
  $G \action X$ is continuous, and $\calS$ is the family of sets of
  the form $\union[n \in \N][\inverse{S_n} S_n \setminus \set{\identity
  {G}}]$, where $S_n \subseteq G$ has cardinality $n$ for all $n
  \in \N$. Then the following are equivalent:
  \begin{enumerate}
    \item There exists $S \in \calS$ for which the family $\principleF
      {S}$ is not in the recurrence spectrum of $G \action X$.
    \item There exist sets $S_n \subseteq G$ of cardinality $n$ and a
      smooth $\aleph_0$-universally \Baire superequivalence relation
      $F$ of $\orbitequivalencerelation{G}{X}$ for which each action
      $G \action \equivalenceclass{x}{F}$ has a $\sigma$-complete
      $\aleph_0$-universally-\Baire set that is $S_n$-wandering for all
      $n \in \N$.
    \item The action $G \action X$ has a very-weakly-wandering
      $\sigma$-complete $\Sigmaclass[0][2]$ set.  
  \end{enumerate}
\end{proposition}

\begin{propositionproof}
  Observe that if $S_n \subseteq G$ for all $n \in \N$, then a set
  $Y \subseteq X$ is $S_n$-wandering for all $n \in \N$ if and
  only if it is $(\union[n \in \N][\inverse{S_n} S_n \setminus
  \set{\identity{G}}])$-transient, and appeal to Proposition \ref
  {transience:calStransient:global}.
\end{propositionproof}

Although we are already in position to establish the analog of
Proposition \ref{transience:weakwandering:local} for very weak
wandering, the following observations will allow us to obtain a
substantially stronger complexity bound.

\begin{proposition} \label{transience:openneighborhoods}
  Suppose that $G$ is a topological group, $X$ is a topological
  space, $G \action X$ is continuous, and $U \subseteq X$ is a
  non-empty open set. Then there exist a non-empty open set $V
  \subseteq U$ and an open neighborhood $W \subseteq G$ of
  $\identity{G}$ for which $W \returntimes{V} \inverse{W} \subseteq
  \returntimes{U}$.
\end{proposition}

\begin{propositionproof}
  The continuity of $G \action X$ yields a non-empty open set $V
  \subseteq U$ and an open neighborhood $W \subseteq G$ of
  $\identity{G}$ for which $WV \subseteq U$. To see that $W
  \returntimes{V} \inverse{W} \subseteq U$, note that if $g \in
  \returntimes{V}$ and $g_x, g_y \in W$, then there exist $x, y \in V$
  for which $g \cdot x = y$, in which case the points $x' = g_x \cdot
  x$ and $y' = g_y \cdot y$ are in $U$, so the fact that $g_y g
  \inverse{g_x} \cdot x' = y'$ ensures that $g_y g \inverse{g_x} \in
  \returntimes{U}$.
\end{propositionproof}

\begin{proposition} \label{transience:dense}
  Suppose that $G$ is a topological group, $X$ is a topological
  space, $G \action X$ is continuous, and $S \subseteq G$. Then
  every $S$-wandering non-empty open set $U \subseteq X$ has a
  non-empty open subset $V \subseteq U$ such that for all dense sets
  $H \subseteq G$, there is an injection $\phi \from S \to H$ with the
  property that $V$ is $\image{\phi}{S}$-wandering.
\end{proposition}

\begin{propositionproof}
  By Proposition \ref{transience:openneighborhoods}, there exist a
  non-empty open set $V \subseteq U$ and an open neighborhood
  $W \subseteq G$ of $\identity{G}$ for which $W \returntimes{V} \inverse{W}
  \subseteq \returntimes{U}$. Note that if $g,h \in S$ and $\inverse
  {(gW)}(hW) \intersection \returntimes{V} \neq \emptyset$, then
  the fact that $\inverse{(gW)}(hW) = \inverse{W} \inverse{g} h W$
  yields that $\inverse{g} h \in W \returntimes{V} \inverse{W}
  \subseteq \returntimes{U}$, thus $g = h$. But if $H \subseteq G$
  is dense, then there is a function $\phi \from S \to H$ with
  the property that $\phi(g) \in gW$ for all $g \in S$, and it follows
  that $\phi$ is injective and $V$ is $\image{\phi}{S}$-wandering.
\end{propositionproof}

Given a superequivalence relation $E$ of $\orbitequivalencerelation
{G}{X}$, we say that a set $Y \subseteq X$ is \definedterm{$E$-locally
very-weakly-wandering} if its intersection with each $E$-class is
very weakly wandering.

\begin{proposition} \label{transience:veryweakwandering:local}
  Suppose that $G$ is a \Polish group, $X$ is a \Polish space, $G
  \action X$ is continuous, and $\calS$ is the family of sets of the
  form $\union[n \in \N][\inverse{S_n} S_n \setminus \set{\identity
  {G}}]$, where $S_n \subseteq G$ has cardinality $n$ for all $n
  \in \N$. Then the following are equivalent:
  \begin{enumerate}
    \item The family $\principleF{\calS}$ is not in the recurrence
      spectrum of $G \action X$.
    \item There is a smooth $\aleph_0$-universally \Baire
      superequivalence $F$ of $\orbitequivalencerelation{G}{X}$ for
      which each action $G \action \equivalenceclass{x}{F}$ has an
      $\orbitequivalencerelation{G}{X}$-locally very-weakly-wandering
      $\sigma$-complete $\aleph_0$-universally-\Baire set.
    \item The action $G \action X$ has an $\minimaldecomposition{G}
      {X}$-locally very-weakly-wandering $\sigma$-complete
      $\Sigmaclass[0][4]$ set.
  \end{enumerate}
\end{proposition}

\begin{propositionproof}
  As a set $Y \subseteq X$ is $\orbitequivalencerelation{G}
  {X}$-locally very-weakly-wandering if and only if $\returntimes{C
  \intersection Y} \notin \principleF{\calS}$ for all equivalence classes
  $C$ of $\orbitequivalencerelation{G}{X}$, Proposition \ref
  {recurrence:characterization:local} yields $(1) \iff (2)$. The fact that
  every $\minimaldecomposition{G}{X}$-locally
  very-weakly-wandering set $Y \subseteq X$ is
  $\orbitequivalencerelation{G}{X}$-locally very-weakly-wandering
  yields $(3) \implies (2)$. To see $(1) \implies (3)$, fix a countable
  dense set $H \subseteq G$, and let $\calT$ denote the family of
  sets of the form $\union[n \in \N][\inverse{T_n} T_n \setminus \set
  {\identity{H}}]$, where $T_n \subseteq H$ has cardinality $n$ for
  all $n \in \N$. Now observe that if condition (1) holds, then
  Proposition \ref{recurrence:characterization:global} ensures that
  there is no equivalence class $C$ of $\minimaldecomposition{G}
  {X}$ for which $G \action C$ is $\principleF{\calS}$-recurrent, so
  Proposition \ref{transience:dense} implies that there is no
  equivalence class $C$ of $\minimaldecomposition{G}{X}$ for which
  $G \action C$ is $\principleF{\calT}$-recurrent, thus $\principleF
  {\calT}$ is not in the recurrence spectrum of $G \action X$. Note
  that if $U \subseteq G \times X$ and $x \in X$, then
  $\horizontalsection{U}{x} \in \principleF{\calT} \iff \exists n \in \N
  \forall \sequence{h_i}[i < n] \in \functions{n}{H} \exists i \neq j \ (h_i
  = h_j \mathor \inverse{h_i} h_j \in \horizontalsection{U}{x})$, so
  $\principleF{\calT}$ is $\Sigmaclass[0][3]$-on-open. As every set
  $Y \subseteq X$ with the property that $\returntimes{C \intersection
  Y} \notin \principleF{\calT}$ for all equivalence classes $C$ of
  $\minimaldecomposition{G}{X}$ is $\minimaldecomposition{G}
  {X}$-locally very-weakly-wandering, Proposition \ref
  {recurrence:complete} yields condition (3).
\end{propositionproof}

We say that a set $Y \subseteq X$ is \definedterm{non-trivially
transient} if there is a non-pre-compact set $S \subseteq G$
for which $Y$ is $S$-transient.

\begin{proposition} \label{transience:minimallytransient:global}
  Suppose that $G$ is a separable group, $X$ is a \Polish
  space, $G \action X$ is continuous, and $\calS$ is the family
  of non-pre-compact sets $S \subseteq G$. Then the following
  are equivalent:
  \begin{enumerate}
    \item There exists $S \in \calS$ for which the family $\principleF
      {S}$ is not in the recurrence spectrum of $G \action X$.
    \item There exist a non-pre-compact set $S \subseteq G$ and
      a smooth $\aleph_0$-universally \Baire superequivalence $F$
      of $\orbitequivalencerelation{G}{X}$ for which each action $G
      \action \equivalenceclass{x}{F}$ has an $S$-transient
      $\sigma$-complete $\aleph_0$-universally-\Baire set.
    \item The action $G \action X$ has a non-trivially-transient
      $\sigma$-complete $\Sigmaclass[0][2]$ set.
  \end{enumerate}
\end{proposition}

\begin{propositionproof}
  Observe that a set $Y \subseteq X$ is non-trivially transient if and
  only if it is $\calS$-transient, and appeal to Proposition \ref
  {transience:calStransient:global}.
\end{propositionproof}

Although we are already in position to establish the analog of
Proposition \ref{transience:weakwandering:local} for non-trivial
transience, we will again first establish an observation yielding
a substantially stronger complexity bound.

\begin{proposition} \label{transience:nonprecompact}
  Suppose that $G$ is a locally compact group, $X$ is a
  topological space, $G \action X$ is continuous, and $S \subseteq
  G$ is not pre-compact. Then every $S$-transient non-empty open
  set $U \subseteq X$ has a non-empty open subset $V \subseteq
  U$ such that for all dense sets $H \subseteq G$, there is a function
  $\phi \from S \to H$ with the property that $\image{\phi}{S}$ is
  not pre-compact and $V$ is $\image{\phi}{S}$-transient.
\end{proposition}

\begin{propositionproof}
  By Proposition \ref{transience:openneighborhoods}, there exist a
  non-empty open set $V \subseteq U$ and an open neighborhood
  $W \subseteq G$ of $\identity{G}$ for which $\returntimes{V}
  \inverse{W} \subseteq \returntimes{U}$. As $G$ is locally compact,
  we can assume that $W$ is pre-compact. As $H \subseteq G$ is
  dense, there is a function $\phi \from S \to H$ with the property
  that $\phi(g) \in gW$ for all $g \in S$. Then $S \subseteq \image
  {\phi}{S} \inverse{W}$, so $\image{\phi}{S}$ is not pre-compact.
  And $\returntimes{V} \inverse{W} \intersection S \subseteq
  \returntimes{U} \intersection S = \emptyset$, so $\returntimes{V}
  \intersection \image{\phi}{S} \subseteq \returntimes{V} \intersection
  SW = \emptyset$, thus $V$ is $\image{\phi}{S}$-transient.
\end{propositionproof}

Given a superequivalence relation $E$ of $\orbitequivalencerelation
{G}{X}$, we say that a set $Y \subseteq X$ is \definedterm
{$E$-locally non-trivially-transient} if its intersection with each
$E$-class is non-trivially transient.

\begin{proposition} \label{transience:minimallytransient:local}
  Suppose that $G$ is a locally compact \Polish group, $X$ is a
  \Polish space, $G \action X$ is continuous, $\calF$ is the family of
  co-pre-compact sets $S \subseteq G$ containing $\identity{G}$,
  and $\calS$ is the family of non-pre-compact sets $S \subseteq G$.
  Then the following are equivalent:
  \begin{enumerate}
    \item The family $\calF$ is not in the recurrence spectrum of
      $G \action X$.
    \item The family $\principleF{\calS}$ is not in the recurrence
      spectrum of $G \action X$.
    \item There is a smooth $\aleph_0$-universally \Baire
      superequivalence $F$ of $\orbitequivalencerelation{G}{X}$ for
      which each action $G \action \equivalenceclass{x}{F}$ has an
      $\orbitequivalencerelation{G}{X}$-locally non-trivially-transient
      $\sigma$-complete $\aleph_0$-universally-\Baire set.
    \item The action $G \action X$ has an $\minimaldecomposition{G}
      {X}$-locally non-trivially-transient $\sigma$-complete $\Sigmaclass
      [0][4]$ set.
  \end{enumerate}
\end{proposition}

\begin{propositionproof}
  As $\principleF{\calS}$ is the family of co-pre-compact sets $S
  \subseteq G$, the fact that $\identity{G} \in \returntimes{Y}$ for all
  $Y \subseteq X$ yields $(1) \iff (2)$. As a set $Y \subseteq X$
  is $\orbitequivalencerelation{G}{X}$-locally non-trivially-transient if
  and only if $\returntimes{C \intersection Y} \notin \principleF{\calS}$
  for all equivalence classes $C$ of $\orbitequivalencerelation{G}
  {X}$, Proposition \ref{recurrence:characterization:local} yields $(2)
  \iff (3)$. The fact that every $\minimaldecomposition{G}{X}$-locally
  non-trivially-transient set $Y \subseteq X$ is
  $\orbitequivalencerelation{G}{X}$-locally non-trivially-transient
  yields $(4) \implies (2)$. To see $(2) \implies (4)$, fix a countable
  dense set $H \subseteq G$, and let $\calT$ denote the family of
  non-pre-compact sets $T \subseteq G$ that are moreover
  contained in $H$. Now observe that if condition (2) holds, then
  Proposition \ref{recurrence:characterization:global} ensures that
  there is no equivalence class $C$ of $\minimaldecomposition{G}
  {X}$ for which $G \action C$ is $\principleF{\calS}$-recurrent, so
  Proposition \ref{transience:nonprecompact} implies that there is no
  equivalence class $C$ of $\minimaldecomposition{G}{X}$ for
  which $G \action C$ is $\principleF{\calT}$-recurrent, thus
  $\principleF{\calT}$ is not in the recurrence spectrum of $G \action
  X$. Fix an increasing sequence $\sequence{K_n}[n \in
  \N]$ of compact subsets of $G$ that is \definedterm{cofinal} in the
  sense that every compact set $K \subseteq G$ is contained in
  some $K_n$, and note that if $U \subseteq G \times X$ and $x \in
  X$, then $\horizontalsection{U}{x} \in \principleF{\calT} \iff \exists
  n \in \N \ H \subseteq K_n \union \horizontalsection{U}{x}$, so
  $\principleF{\calT}$ is $\Sigmaclass[0][3]$-on-open. As every set
  $Y \subseteq X$ with the property that $\returntimes{C \intersection
  Y} \notin \principleF{\calT}$ for all equivalence classes $C$ of
  $\minimaldecomposition{G}{X}$ is $\minimaldecomposition{G}
  {X}$-locally non-trivially-transient, Proposition \ref
  {recurrence:complete} yields condition (4).
\end{propositionproof}

\section{Generic compressibility}
\label{measures}

Suppose that $E$ is a \Borel equivalence relation on $X$ that is
\definedterm{countable}, in the sense that all of its equivalence
classes are countable. We say that a function $\rho \from E \to
\openinterval{0}{\infty}$ is a \definedterm{cocycle} if $\rho(x, z) =
\rho(x, y) \rho(y, z)$ whenever $x \mathrel{E} y \mathrel{E} z$. When
$\rho \from E \to \openinterval{0}{\infty}$ is a \Borel cocycle, we say
that a \Borel measure $\mu$ on $X$ is \definedterm{$\rho$-invariant}
if $\mu(\image{T}{B}) = \int_B \rho(T(x), x) \ d\mu(x)$ for all \Borel
sets $B \subseteq X$ and \Borel automorphisms $T \from X \to X$
such that $\graph{T} \subseteq E$. We say that $\rho$ is
\definedterm{aperiodic} if $\sum_{y \in \equivalenceclass{x}{E}} \rho
(y, x) = \infty$ for all $x \in X$. Here we generalize the following fact
to orbit equivalence relations induced by \Borel actions of locally
compact \Polish groups, while simultaneously strengthening
comeagerness to a notion under which the recurrence spectrum is
invariant.

\begin{theorem}[\Kechris-\Miller] \label{measures:KM04}
  Suppose that $X$ is a standard \Borel space, $E$ is a countable
  \Borel equivalence relation on $X$, and $\rho \from E \to
  \openinterval{0}{\infty}$ is an aperiodic \Borel cocycle. Then there
  is an $E$-invariant comeager \Borel set $C \subseteq X$ that is
  null with respect to every $\rho$-invariant \Borel probability measure.
\end{theorem}

A function $\phi \from X \to Z$ is \definedterm{$E$-invariant} if $\phi
(x) = \phi(y)$ whenever $x \mathrel{E} y$. The \definedterm
{$E$-saturation} of a set $Y \subseteq X$ is the set of $x \in X$ for
which there exists $y \in Y$ such that $x \mathrel{E} y$. We say that
a \Borel probability measure $\mu$ on $X$ is \definedterm
{$E$-quasi-invariant} if the $E$-saturation of every $\mu$-null set
$N \subseteq X$ is $\mu$-null. Let $\probabilitymeasures{X}$
denote the standard \Borel space of \Borel probability measures
on $X$ (see, for example, \cite[\S17.E]{Kec95}). The \definedterm
{push-forward} of a \Borel measure $\mu$ on $X$ through a \Borel
function $\phi \from X \to Y$ is the \Borel measure $\pushforward
{\phi}{\mu}$ on $Y$ given by $(\pushforward{\phi}{\mu})(B) = \mu
(\preimage{\phi}{B})$ for all \Borel sets $B \subseteq Y$.

\begin{proposition} \label{measures:cocycle}
  Suppose that $X$ is a standard \Borel space, $E$ is a countable
  \Borel equivalence relation on $X$, and $\phi \from X \to
  \probabilitymeasures{X}$ is an $E$-invariant \Borel function such
  that $\mu$ is $E$-quasi-invariant and $\preimage{\phi}{\mu}$ is
  $\mu$-conull for all $\mu \in \image{\phi}{X}$. Then there is a
  \Borel cocycle $\rho \from E \to \openinterval{0}{\infty}$ such that
  $\mu$ is $\rho$-invariant for all $\mu \in \image{\phi}{X}$.
\end{proposition}

\begin{propositionproof}
  By standard change of topology results (see, for example, \cite
  [\S13]{Kec95}), we can assume that $X$ is a zero-dimensional
  \Polish space. Fix a compatible complete ultrametric on $X$. By
  \cite[Theorem 1]{FM77}, there is a countable group
  $\Gamma$ of \Borel automorphisms of $X$ whose induced
  orbit equivalence relation is $E$. For all $\gamma \in \Gamma$,
  define $\rho_\gamma \from X_\gamma \to \openinterval{0}{\infty}$
  by $\rho_\gamma(x) = \lim_{\epsilon \goesto 0} (\pushforward
  {(\inverse{\gamma})}{\phi(x)})(\ball{x}{\epsilon}) / \phi(x)(\ball{x}
  {\epsilon})$, where $X_\gamma$ is the set of $x \in X$ for which
  this limit exists and lies in $\openinterval{0}{\infty}$.
  
  Note that if $\gamma \in \Gamma$, $\mu \in \image{\phi}{X}$,
  and $\psi \from X \to \openinterval{0}{\infty}$ is a \Radon-\Nikodym
  derivative of $\pushforward{(\inverse{\gamma})}{\mu}$ with respect
  to $\mu$ (see, for example, \cite[\S17.A]{Kec95}), then the
  straightforward generalization of the \Lebesgue density theorem
  for \Polish ultrametric spaces (see, for example, \cite[Proposition
  2.10]{Mil08a}) to integrable functions ensures that $\psi(x)
  = \lim_{\epsilon \goesto 0} \int_{\ball{x}{\epsilon}} \psi \ d\mu / \mu
  (\ball{x}{\epsilon}) = \rho_\gamma(x)$ for $\mu$-almost all $x \in X$.
  
  It immediately follows that for all $\gamma \in \Gamma$, the
  complement of $X_\gamma$ is null with respect to every $\mu
  \in \image{\phi}{X}$. Moreover, if $B \subseteq X$ is \Borel,
  $\gamma, \delta \in \Gamma$, and $\mu \in \image{\phi}{X}$, then
  \begin{align*}
    \pushforward{\inverse{(\gamma \delta)}}{\mu}(B)
      & = \int_{\image{\delta}{B}} \rho_\gamma(x) \ d\mu(x) \\
      & = \int_B \rho_\gamma(\delta \cdot x) \ d(\pushforward{(\inverse
        {\delta})}{\mu})(x) \\
      & = \int_B \rho_\gamma(\delta \cdot x) \rho_\delta(x) \ d\mu(x),
  \end{align*}
  so the almost-everywhere uniqueness of \Radon-\Nikodym
  derivatives ensures that the set of $x \in X$ for which there exist
  $\gamma, \delta \in \Gamma$ such that $\rho_{\gamma \delta}(x)
  \neq \rho_\gamma(\delta \cdot x) \rho_\delta(x)$ is null with respect
  to every $\mu \in \image{\phi}{X}$.
  
  Let $N$ denote the $E$-saturation of the union of these sets, and
  let $\rho \from E \to \openinterval{0}{\infty}$ be the extension of the
  constant cocycle on $\restriction{E}{N}$ given by $\rho(\gamma
  \cdot x, x) = \rho_\gamma(x)$ for all $\gamma \in \Gamma$ and $x
  \in X$.
\end{propositionproof}

As a consequence, we obtain the following.

\begin{theorem} \label{measures:comeagernullset}
  Suppose that $X$ is a \Polish space, $E$ is a \Borel equivalence
  relation on $X$ admitting a \Borel complete set $B \subseteq X$
  on which $E$ is countable, $F$ is a superequivalence relation of
  $E$ for which every $F$-class is \Gdelta and the $F$-saturation
  of every open set is \Borel, and $\phi \from X \to \probabilitymeasures
  {X}$ is an $E$-invariant \Borel function for which every measure
  $\mu \in \image{\phi}{X}$ has $\mu$-conull $\phi$-preimage and 
  concentrates off of \Borel sets on which $E$ is smooth. Then there
  is an $E$-invariant \Borel set $C \subseteq X$ that is comeager in
  every $F$-class, but null with respect to every measure in $\image
  {\phi}{X}$.
\end{theorem}

\begin{theoremproof}
  By the \Lusin-\Novikov uniformization theorem (see, for example,
  \cite[Theorem 18.10]{Kec95}), there is a \Borel extension $\pi
  \from X \to B$ of the identity function on $B$ whose graph is
  contained in $E$. Fix a sequence $\sequence{\epsilon_n}[n \in
  \N]$ of positive real numbers whose sum is $1$, in addition to a
  countable group $\set{\gamma_n}[n \in \N]$ of \Borel
  automorphisms of $B$ whose induced orbit equivalence relation is
  $\restriction{E}{B}$, and define $\psi \from B \to
  \probabilitymeasures{B}$ by $\psi(x) = \sum_{n \in \N} \pushforward
  {(\gamma_n \composition \pi)}{\phi(x)} / \epsilon_n$. As each $\nu
  \in \image{\psi}{B}$ is $(\restriction{E}{B})$-quasi-invariant,
  Proposition \ref{measures:cocycle} yields a \Borel cocycle $\rho
  \from \restriction{E}{B} \to \openinterval{0}{\infty}$ such that every
  $\nu \in \image{\psi}{B}$ is $\rho$-invariant.
  
  Given $\nu \in \image{\psi}{B}$, fix $x \in B$ for which $\nu =
  \psi(x)$, set $\mu = \phi(x)$, and observe that $\nu(\preimage{\psi}{\nu})
  \ge \mu(\preimage{\phi}{\mu}) = 1$. Moreover, as $E$ is smooth on
  the periodic part $P = \set{x \in B}[\sum_{y \in \equivalenceclass
  {x}{\restriction{E}{B}}} \rho(y, x) < \infty]$ of $\rho$ (see, for example,
  \cite[Proposition 2.1.1]{Mil08b}), and therefore on
  its $E$-saturation, it follows that $\saturation{P}{E}$ is null with
  respect to every measure in $\image{\phi}{X}$, thus $P$ is null with
  respect to every measure in $\image{\psi}{B}$.
  
  By the proof of Theorem \ref{measures:KM04} (see \cite
  [Theorem 13.1]{KM04}), there is a \Borel set $R \subseteq
  \Bairespace \times B$, whose vertical sections are $(\restriction{E}
  {B})$-invariant and null with respect to every $\rho$-invariant \Borel
  probability measure, such that every $x \in B$ is contained
  in comeagerly-many vertical sections of $R$. It follows that the
  vertical sections of the set $S = \preimage{(\id \times \pi)}{R}$ are
  $E$-invariant and null with respect to every measure in $\image
  {\phi}{X}$, and every $x \in X$ is contained in comeagerly-many
  vertical sections of $S$. The \Kuratowski-\Ulam theorem therefore
  ensures that for all $x \in X$, comeagerly-many vertical sections of
  $S$ are comeager in $\equivalenceclass{x}{F}$. 
  
  By \cite{Sri79}, there is a \Borel set $D \subseteq X$
  intersecting every $F$-class in a single point. As the $F$-saturation
  of every open set is \Borel, the usual proof of the \Montgomery-\Novikov
  theorem that the pointclass of \Borel sets is closed under category
  quantifiers (see, for example, \cite[Theorem 16.1]{Kec95}) shows that
  $\set{\pair{b}{x} \in \Bairespace \times X}[\verticalsection{S}{b} \text{ is
  comeager in } \equivalenceclass{x}{F}]$ and $\set{\pair{b}{x} \in
  \Bairespace \times X}[\verticalsection{S}{b} \text{ is non-meager in }
  \equivalenceclass{x}{F}]$ are \Borel, so \cite[Theorem 18.6]{Kec95}
  yields a \Borel function $\beta \from D \to \Bairespace$ such that
  $\verticalsection{S}{\beta(x)}$ is comeager in $\equivalenceclass
  {x}{F}$ for all $x \in D$. Then the set $C = \union[x \in D][\verticalsection
  {S}{\beta(x)} \intersection \equivalenceclass{x}{F}]$ is as desired.
\end{theoremproof}

We say that a function $\rho \from G \times X \to \openinterval{0}
{\infty}$ is a \definedterm{cocycle} if $\rho(gh, x) = \rho(g, h \cdot x)
\rho(h, x)$ for all $g, h \in G$ and $x \in X$. When $\rho \from G
\times X \to \openinterval{0}{\infty}$ is a \Borel cocycle, we say that
a \Borel measure $\mu$ on $X$ is \definedterm{$\rho$-invariant} if
$\mu(gB) = \int_B \rho(g, x) \ d\mu(x)$ for all \Borel sets $B
\subseteq X$ and group elements $g \in G$. The following fact is the
desired generalization of Theorem \ref{measures:KM04}.

\begin{theorem} \label{measures:genericcompressibility}
  Suppose that $G$ is a locally compact \Polish group, $X$ is a \Polish
  space, $G \action X$ is a continuous action, $F$ is a superequivalence
  relation of $\orbitequivalencerelation{G}{X}$ for which every $F$-class
  is \Gdelta and the $F$-saturation of every open set is \Borel, and $\rho
  \from G \times X \to \openinterval{0}{\infty}$ is a \Borel cocycle with
  the property that every $G$-orbit is null with respect to every
  $\rho$-invariant \Borel probability measure. Then there is a 
  $G$-invariant \Borel set $C \subseteq X$ that is comeager in every
  $F$-class, but null with respect to
  every $\rho$-invariant \Borel probability measure.
\end{theorem}

\begin{theoremproof}
  By \cite[Theorem 1.1]{Kec92}, there is a complete
  \Borel set $B \subseteq X$ on which $\orbitequivalencerelation{G}{X}$
  is countable. Fix a $\rho$-invariant uniform ergodic decomposition
  $\phi \from X \to \probabilitymeasures{X}$ of $G \action X$ (see \cite
  [Theorem 5.2]{GS00}), and appeal to Theorem \ref
  {measures:comeagernullset}.
\end{theoremproof}

We next check that the special case of Theorem \ref
{measures:genericcompressibility} for $F = \minimaldecomposition
{G}{X}$ provides a proper strengthening of Theorem \ref
{measures:KM04}. While this can be seen as a consequence
of the \Kuratowski-Ulam theorem, we will show that the usual proof
of the latter easily adapts to yield a generalization to a natural class of
equivalence relations containing $\minimaldecomposition{G}{X}$.

\begin{theorem}
  Suppose that $X$ is a second-countable \Baire space, $E$ is an
  equivalence relation on $X$ such that every $E$-class is a \Baire
  space and the $E$-saturation of every open subset of $X$ is open,
  and $B \subseteq X$ has the \Baire property. Then:
  \begin{enumerate}
    \item $\forcomeagerlymany x \in X \ B$ has the \Baire property in
      $\equivalenceclass{x}{E}$.
    \item $B$ is comeager $\iff \forcomeagerlymany x \in X \ B$ is
      comeager in $\equivalenceclass{x}{E}$.
  \end{enumerate}
\end{theorem}

\begin{theoremproof}
  We begin with a simple observation.
  
  \begin{lemma} \label{measures:densesaturation}
    Suppose that $U \subseteq X$ is a non-empty open set and
    $V \subseteq U$ is a dense open set. Then $\saturation{V}{E}$
    is dense in $\saturation{U}{E}$.
  \end{lemma}
  
  \begin{lemmaproof}
    If $W \subseteq X$ is open and $\saturation{V}{E} \intersection W
    = \emptyset$, then $V \intersection \saturation{W}{E} = \emptyset$,
    so the openness of $\saturation{W}{E}$ ensures that $\closure{V}
    \intersection \saturation{W}{E} = \emptyset$, thus the density of
    $V$ implies that $U \intersection \saturation{W}{E} = \emptyset$,
    hence $\saturation{U}{E} \intersection W = \emptyset$.
  \end{lemmaproof}
  
  To see the special case of ($\implies$) of (2) when $B \subseteq X$
  is open, note that if $U \subseteq X$ is non-empty and open, then
  Lemma \ref{measures:densesaturation} yields that $\saturation{B
  \intersection U}{E}$ is dense in $\saturation{U}{E}$, and therefore
  $\forcomeagerlymany x \in X \ (x \in \saturation{U}{E} \implies x \in
  \saturation{B \intersection U}{E})$, or equivalently,
  $\forcomeagerlymany x \in X \ (U \intersection \equivalenceclass{x}{E}
  \neq \emptyset \implies B \intersection U \intersection
  \equivalenceclass{x}{E} \neq \emptyset)$. As $X$ is second countable,
  it follows that $\forcomeagerlymany x \in X \ B$ is dense in
  $\equivalenceclass{x}{E}$.
  
  To see ($\implies$) of (2), suppose that $B \subseteq X$
  is comeager, fix dense open sets $B_n \subseteq X$ for which
  $\intersection[n \in \N][B_n] \subseteq B$, and appeal to the
  special case for open sets to obtain that $\forcomeagerlymany x \in
  X \ \intersection[n \in \N][B_n]$ is comeager in $\equivalenceclass
  {x}{E}$.
  
  To see (1), fix an open set $U \subseteq X$ for which $B
  \symmetricdifference U$ is meager, and note that
  $\forcomeagerlymany x \in X \ B \symmetricdifference U$ is meager
  in $\equivalenceclass{x}{E}$, by $(\implies)$ of (2).
  
  To see ($\impliedby$) of (2), suppose that $B$ is not
  comeager, fix a non-empty open set $V \subseteq X$ in which
  $B$ is meager, note that $\forall x \in V \ V \intersection
  \equivalenceclass{x}{E} \neq \emptyset$, and appeal to
  ($\implies$) of (2) to obtain that $\forcomeagerlymany
  x \in X \ B \intersection V$ is meager in $\equivalenceclass{x}{E}$,
  thus $\forcomeagerlymany x \in V \ B$ is not comeager in
  $\equivalenceclass{x}{E}$.
\end{theoremproof}

Finally, we check that no condition on the recurrence spectrum
can yield the existence of an invariant \Borel probability measure.

\begin{theorem} \label{recurrence:compressibility}
  Suppose that $G$ is a locally compact \Polish group, $X$ is a
  standard \Borel space, $G \action X$ is \Borel, and $\rho \from
  G \times X \to \openinterval{0}{\infty}$ is a \Borel cocycle for which
  every $G$-orbit is null with respect to every $\rho$-invariant \Borel
  probability measure. Then there is a $G$-invariant \Borel set
  $B \subseteq X$ that is null with respect to every $\rho$-invariant
  \Borel probability measure but for which the recurrence spectra of
  $G \action B$ and $G \action X$ coincide.
\end{theorem}

\begin{theoremproof}
  We can assume that $X$ is a \Polish space and $G \action X$ is
  continuous. By Theorem \ref{measures:genericcompressibility},
  there is an $\orbitequivalencerelation{G}{X}$-invariant \Borel set
  $B \subseteq X$ that is comeager in every $\minimaldecomposition
  {G}{X}$-class, but null with respect to every $\rho$-invariant \Borel
  probability measure. Proposition \ref
  {recurrence:characterization:global} then ensures that $G \action
  B$ and $G \action X$ have the same recurrence spectra.
\end{theoremproof}

\subsection*{Acknowledgements}
We would like to thank the anonymous referee for her suggestions.

\bibliographystyle{asl}
\bibliography{bibliography}

\end{document}